\documentclass[a4paper,12pt]{amsart}
\usepackage{amssymb,amsfonts,amsxtra,amscd,mathrsfs}
\usepackage[all]{xy}
\usepackage{fullpage}
\theoremstyle{plain}
\newtheorem{theorem}{Theorem}[section]
\newtheorem{lemma}[theorem]{Lemma}
\newtheorem{cor}[theorem]{Corollary}
\newtheorem{prop}[theorem]{Proposition}
\theoremstyle{definition}
\newtheorem{defi}[theorem]{Definition}
\newtheorem{example}[theorem]{Example}
\theoremstyle{remark}
\newtheorem{rem}[theorem]{Remark}
\numberwithin{equation}{section}
\newcommand{\galg}[1]{\ensuremath{\mathfrak{g}_{#1}}}
\newcommand{\tgalg}[1]{\ensuremath{\tilde{\mathfrak{g}}_{#1}}}
\newcommand{\osp}[1]{\ensuremath{\mathfrak{osp}_{#1}}}
\newcommand{\gc}[1][\bullet\bullet]{\ensuremath{\mathcal{G}_{#1}}}
\newcommand{\gh}[1][\bullet\bullet]{\ensuremath{H\mathcal{G}_{#1}}}
\newcommand{\dilim}[2]{\ensuremath{\varinjlim_{#1} #2}}
\newcommand{\noproof}{\begin{flushright} \ensuremath{\square} \end{flushright}}
\DeclareMathOperator{\sgn}{sgn}
\DeclareMathOperator{\id}{id}
\begin{document}

\begin{abstract}
In this paper we will prove a super-analogue of a well-known result by Kontsevich which states that the homology of a certain complex which is generated by isomorphism classes of oriented graphs can be calculated as the Lie algebra homology of an infinite-dimensional Lie algebra of symplectic vector fields.
\end{abstract}
\title{A super-analogue of Kontsevich's theorem on graph homology}
\author{Alastair Hamilton}
\address{Mathematics Department, Bristol University, Bristol, England. BS8 1TW.} \email{a.e.hamilton@bristol.ac.uk}
\keywords{Graph, Homology, Lie superalgebra, Moduli space, Invariant theory.}
\subjclass[2000]{17B56, 17B65, 17B66, 57M27.}
\thanks{The work of the author was partially supported by an EPSRC grant No. GR/SO7148/01.}
\maketitle

\section{Introduction}

Graph homology is a subject which is at the junction of several areas of mathematics including: invariants of differentiable manifolds, infinity-algebras, algebraic geometry and quantum field theory. There is a version of the graph complex for every cyclic operad \cite{vogtmann}; the complex is generated by isomorphism classes of decorated oriented graphs and the differential in this complex sends such a graph to an alternating sum of graphs which are obtained by contracting the edges of the initial graph.

In this paper we shall only consider the graph complex associated to the commutative operad. This is done for the sake of simplicity and brevity; however the treatment that we give here goes through in other important cases with minor alterations. In this context a `commutative' graph is the same thing as a one-dimensional cell complex. The corresponding graph complex is related to knot invariants and invariants of differentiable manifolds; cf. \cite{kontfd}, \cite{schwarz} and also see \cite{bott}. Graph complexes corresponding to other cyclic operads are also related to some interesting moduli spaces; for instance the graph complex associated to the associative operad is related to moduli spaces of complex algebraic curves, cf. \cite{kontfd} and \cite{penner}.

The objective of this paper is to prove a super-analogue of the famous theorem by Kontsevich \cite{kontsg} which says that the homology of the graph complex can be computed as the stable Lie algebra homology of certain Lie algebras consisting of symplectic vector fields. The super-analogue which we intend to prove is a generalisation of Kontsevich's original result which he alluded to in \cite{kontfd}. The main tool that we need in order to prove it is the invariant theory for the `orthosymplectic' Lie superalgebras $\osp{2n|m}$ which is provided by \cite{sergeev}. A standard modern reference on supermathematics is Deligne and Morgan's article \cite{delsym}.

As well as generalising Kontsevich's original theorem, the treatment that we shall present here has another advantage; it is fairly brief and the calculations are relatively simple. Recently, there has been work done by other authors \cite{vogtmann}, \cite{mahajan} in understanding and generalising Kontsevich's work, however these treatments are quite long and the calculations can be difficult to follow. By working with the \emph{coinvariants} of the Lie superalgebra $\osp{2n|m}$ rather than its \emph{invariants}, it is possible to simplify some of the calculations that are involved.

The main reason however for desiring a super-analogue of Kontsevich's original result concerns the constructions introduced by Kontsevich in \cite{kontfd}. These take infinity-algebras with an invariant inner product and produce (co)homology classes in the graph complex. These constructions have strong links with two aspects of quantum field theory: the Feynman calculus and the Batalin-Vilkovisky formalism. These constructions and their relationship with quantum field theory form the basis of papers which are currently in preparation by the author in collaboration with Andrey Lazarev.

In this paper we choose to work over the field of complex numbers $\mathbb{C}$. This is purely for convenience and the results of this paper hold over any field of characteristic zero, or even any ring containing the field of rational numbers.

The layout of the paper is as follows: in section \ref{sec_gphhom} we recall the basic definition of graph homology and layout the terminology that we will use for the remainder of the paper. Section \ref{sec_liehom} recalls the basic definition of Lie superalgebra homology and introduces a series of Lie superalgebras which form the natural generalisation of the series of Lie algebras introduced by Kontsevich in \cite{kontsg}. Section \ref{sec_invthy} formulates the results on the invariant theory of the Lie superalgebras $\osp{2n|m}$ that we will need for the main theorem, which is proved in section \ref{sec_gphhommain}.

\textit{Acknowledgement}: The author would like to thank Andrey Lazarev for his advice and assistance during the completion of this work.

\subsection{Notation and conventions}

In this paper we will work over the field of complex numbers $\mathbb{C}$. We shall denote the canonical $2n|m$ dimensional vector space by $\mathbb{C}^{2n|m}$ which we shall write in the canonical coordinates familiar to symplectic geometry;
\[ \mathbb{C}^{2n|m} = \langle p_1,\ldots,p_n;q_1,\ldots,q_n;x_1,\ldots,x_m\rangle, \]
where the $p_i$'s and $q_i$'s are even and the $x_i$'s are odd.

The canonical symplectic form on $\mathbb{C}^{2n|m}$ will be denoted simply by $\langle -,- \rangle$; it is given by the formula
\begin{equation} \label{eqn_cansymfrm}
\begin{split}
\langle p_i,q_j \rangle = \langle x_i,x_j \rangle & = \delta_{ij}, \\
\langle x_i,p_j \rangle = \langle x_i,q_j \rangle = \langle p_i,p_j \rangle = \langle q_i,q_j \rangle & = 0.
\end{split}
\end{equation}
The Lie subalgebra of $\mathfrak{gl}_{2n|m}(\mathbb{C})$ consisting of all linear endomorphisms $\xi$ which satisfy the relation
\[ \langle \xi(a), b \rangle + (-1)^{|\xi||a|}\langle a , \xi(b) \rangle = 0;  \text{ for all } a,b \in \mathbb{C}^{2n|m}; \]
will be denoted by $\osp{2n|m}$.

The subgroup of $Gl_{2n|m}(\mathbb{C})$ consisting of all invertible linear endomorphisms $\phi$ satisfying the relation
\[ \langle \phi(a),\phi(b) \rangle = \langle a,b \rangle; \text{ for all } a,b \in \mathbb{C}^{2n|m}; \]
will be denoted by $OSp_{2n|m}$.

Let $\sigma \in S_m$ and $k_1,\ldots,k_m$ be positive integers such that $k_1+\ldots+k_m=N$, we define the permutation $\sigma_{(k_1,\ldots,k_m)} \in S_N$ by the following commutative diagram:
\begin{displaymath}
\xymatrix{ V^{\otimes k_1}\otimes\ldots\otimes V^{\otimes k_m} \ar@{=}[d] \ar^{\sigma}_{x_1\otimes\ldots\otimes x_m \mapsto x_{\sigma(1)}\otimes\ldots\otimes x_{\sigma(m)}}[rrrr] &&&& V^{\otimes k_{\sigma(1)}}\otimes\ldots\otimes V^{\otimes k_{\sigma(m)}} \ar@{=}[d] \\ V^{\otimes N} \ar^{\sigma_{(k_1,\ldots,k_m)}}_{y_1\otimes\ldots\otimes y_N \mapsto y_{\sigma_{(k_1,\ldots,k_m)}[1]}\otimes\ldots\otimes y_{\sigma_{(k_1,\ldots,k_m)}[N]}}[rrrr] &&&& V^{\otimes N} }
\end{displaymath}

Let $k$ be a positive integer. A partition $c$ of $\{1,\ldots 2k\}$ such that every $x \in c$ is a set consisting of precisely two elements will be called a \emph{chord diagram}. The set of all such chord diagrams will be denoted by $\mathscr{C}(k)$. A chord diagram with an ordering (orientation) of each two point set in that chord diagram will be called an \emph{oriented chord diagram}. The set of all oriented chord diagrams will be denoted by $\mathscr{OC}(k)$.

$S_{2k}$ acts on $\mathscr{OC}(k)$ as follows: given an oriented chord diagram $c:=(i_1,j_1),\ldots,(i_k,j_k)$ and a permutation $\sigma \in S_{2k}$,
\[ \sigma\cdot c := (\sigma(i_1),\sigma(j_1)),\ldots,(\sigma(i_k),\sigma(j_k)). \]
$S_{2k}$ acts on $\mathscr{C}(k)$ in a similar fashion.

A vector superspace $V$ is a vector space with a natural decomposition $V=V_0\oplus V_1$ into an even part ($V_0$) and a odd part ($V_1$). Given a vector superspace $V$, we will denote the parity reversion of $V$ by $\Pi V$. The parity reversion $\Pi V$ is defined by the identities
\[ (\Pi V)_0 := V_1, \quad (\Pi V)_1 := V_0. \]
The canonical bijection $V \to \Pi V$ will be denoted by $\Pi$. In this paper we will refer to all vector superspaces simply as `vector spaces'.

Given a vector superspace $V$ we can define a free supercommutative algebra $S(V)$ as the quotient of the free associative algebra $T(V):=\bigoplus_{k=0}^\infty V^{\otimes k}$ by the ideal generated by the relation
\[ x \otimes y = (-1)^{|x||y|} y \otimes x; \quad x,y \in V. \]
This algebra has a natural $\mathbb{Z}$-grading which we will refer to as the grading by \emph{order}; we say an element $x \in S(V)$ has homogeneous \emph{order} $n$ if $x \in S^n(V)$.

\section{Graph homology} \label{sec_gphhom}

In this section we will recall the basic definitions of graph homology (cf. \cite{kontsg} and \cite{kontfd}). In order to define the graph complex it will be necessary to consider \emph{oriented graphs}; these are graphs with an extra structure called an \emph{orientation}. Informally, the graph complex is a complex which is linearly generated by isomorphism classes of oriented graphs and where the differential is given by taking an alternating sum of edge contractions over all edges of a graph. This section will detail all the prerequisite material needed to make this definition explicit and lay down the appropriate terminology.

We will use the following set theoretic definition of a graph (see Igusa \cite{igusa}):

\begin{defi}
Choose a fixed infinite set $\Xi$ which is disjoint from its power set; a graph $\Gamma$ is a finite subset of $\Xi$ (the set of half-edges) together with:
\begin{enumerate}
\item
a partition $V(\Gamma)$ of $\Gamma$ ($V(\Gamma)$ is the set of vertices of $\Gamma$),
\item
a partition $E(\Gamma)$ of $\Gamma$ into sets having cardinality equal to two ($E(\Gamma)$ is the set of edges of $\Gamma$).
\end{enumerate}
We say that a vertex $v \in V$ has valency $n$ if $v$ has cardinality $n$. The elements of $v$ are called the \emph{incident half-edges} of $v$. An edge whose composite half-edges are incident to the same vertex will be called a loop. In this paper we shall assume that all our graphs have vertices of valency $\geq 3$.
\end{defi}

As already mentioned, in order to define graph homology it will be necessary to consider oriented graphs. An orientation on a graph is defined as follows:

\begin{defi} \label{def_orient}
Let $\Gamma$ be a graph with $m$ vertices and $n$ edges. Consider the set of structures $\Omega$ on $\Gamma$ consisting of:
\begin{enumerate}
\item
a linear ordering of the vertices,
\item
a linear ordering of the half-edges in every edge of $\Gamma$ (i.e. an orientation of each edge as with a directed graph).
\end{enumerate}
Then
\[ G:=S_m \times \underbrace{S_2 \times \ldots \times S_2}_{n \text{ copies}} \]
acts freely and transitively on $\Omega$; the first factor changes the ordering of vertices and the subsequent factors flip the orientation of edges. We then consider the subgroup $K$ of $G$ consisting of permutations having total sign one, i.e. $K$ is the kernel of the map
\begin{displaymath}
\begin{array}{ccc}
S_m \times \underbrace{S_2 \times \ldots \times S_2}_{n \text{ copies}} & \to & \{ 1,-1 \} \\
(\sigma,\tau_1,\ldots,\tau_n) & \mapsto & \sgn\sigma \sgn\tau_1 \ldots \sgn\tau_n \\
\end{array}
\end{displaymath}
An orientation on $\Gamma$ is defined to be an element of $\Omega/K$. This means that there are only two possible choices for an orientation on $\Gamma$.
\end{defi}

We will also need to be able to identify when two graphs are equivalent; this is given by the following definition:

\begin{defi}
Let $\Gamma$ and $\Gamma'$ be two graphs, we say $\Gamma$ is isomorphic to $\Gamma'$ if there exists a bijective map $f:\Gamma \to \Gamma'$ between the half-edges of $\Gamma$ and $\Gamma'$ such that:
\begin{enumerate}
\item
the image of a vertex of $\Gamma$ under $f$ is a vertex of $\Gamma'$,
\item
the image of an edge of $\Gamma$ under $f$ is an edge of $\Gamma'$.
\end{enumerate}

If $\Gamma$ and $\Gamma'$ are oriented graphs then we say that they are isomorphic if there is a map $f$ as above which preserves the orientations of $\Gamma$ and $\Gamma'$; given orientations $\omega$ and $\omega'$ of $\Gamma$ and $\Gamma'$ respectively and a map $f$ as above there is a naturally induced orientation $f(\omega)$ of $\Gamma'$ defined in an obvious way, we say that $f$ preserves the orientations of $\Gamma$ and $\Gamma'$ if $f(\omega)=\omega'$.
\end{defi}

In order to define the differential in the graph complex it will be necessary to describe how to contract an edge in a given graph in a way which is consistent with any orientation that there may be on that graph:

\begin{defi}
Let $\Gamma$ be an oriented graph and $e=(h,h')$ be one of its edges, which we assume is not a loop. The half-edges $h$ and $h'$ will be incident to vertices denoted by $v$ and $v'$ respectively. We may assume that the orientation on $\Gamma$ is represented by an ordering of the vertices of the form
\[ (v,v',v_3,\ldots,v_m) \]
and where the edge $e$ is oriented in the order $e=(h,h')$. We define a new oriented graph $\Gamma/e$ as follows:
\begin{enumerate}
\item
the set of half edges comprising $\Gamma/e$ are the half edges of $\Gamma$ minus $h$ and $h'$,
\item
the set of vertices of $\Gamma/e$ is
\[ ([v-{h}]\cup [v'-{h'}],v_3,\ldots,v_m), \]
\item
The set of edges of $\Gamma/e$ is the set of edges of $\Gamma$ minus the edge $e$,
\item
The orientation on $\Gamma/e$ is given by ordering the vertices as above in (ii); the edges of $\Gamma/e$ are oriented in the same way as the edges of $\Gamma$. It is clear that this definition produces a well defined orientation modulo the actions of the relevant permutation groups.
\end{enumerate}
\end{defi}

This definition is illustrated by the following diagram:

\begin{example}
Contracting an edge $e$ in a graph $\Gamma$:
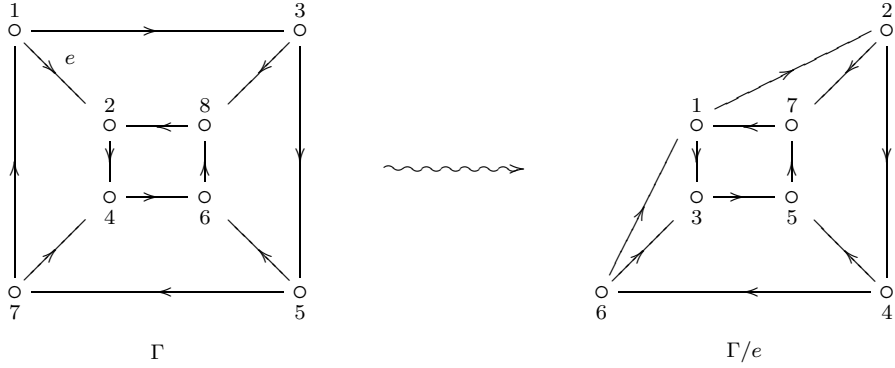
\begin{displaymath}
\xymatrix{
\xygraph{
!~:{@{-}|@{>}}
{\overset{8}{\circ}}="eight":@<+0.9ex>[l]{\overset{2}{\circ}}="two":[d]{\underset{4}{\circ}}="four":@<+0.9ex>[r]{\underset{6}{\circ}}="six":"eight"
[ur]{\overset{3}{\circ}}="three":[ddd]{\underset{5}{\circ}}="five":@<-0.9ex>[lll]{\underset{7}{\circ}}="seven":[uuu]{\overset{1}{\circ}}="one":@<-0.9ex>"three"
"one":@<-0.7ex>^{e}"two" "seven":@<+0.7ex>"four" "five":@<-0.7ex>"six" "three":@<+0.7ex>"eight"
!~:{@{->}} "seven":@{}|{\Gamma}[drrr]
}
&& \ar@<-4ex>@{~>}[rr] &&&&
\xygraph{
!~:{@{-}|@{>}}
{\overset{1}{\circ}}="one":[d]{\underset{3}{\circ}}="three":@<+0.9ex>[r]{\underset{5}{\circ}}="five":[u]{\overset{7}{\circ}}="seven":@<+0.9ex>"one"
:@<-0.3ex>[urr]{\overset{2}{\circ}}="two":[ddd]{\underset{4}{\circ}}="four":@<-0.9ex>[lll]{\underset{6}{\circ}}="six":@<+0.4ex>"one"
"two":@<+0.7ex>"seven" "four":@<-0.6ex>"five" "six":@<+0.6ex>"three"
!~:{@{->}} "six":@{}|{\Gamma/e}[drrr]
}
}
\end{displaymath}
\end{example}

Recall that in Definition \ref{def_orient} it was noted that for any graph $\Gamma$ there were two possible choices of orientations on it, therefore given an orientation $\omega$ of $\Gamma$ we may define $-\omega$ to be the other choice of orientation on $\Gamma$. We are now ready to define the graph complex:

\begin{defi} \label{def_graph}
The underlying space of the graph complex $\gc$ is the free $\mathbb{C}$-vector space generated by isomorphism classes of oriented graphs, modulo the following relation:
\begin{equation} \label{eqn_graphdummy}
(\Gamma,-\omega)=-(\Gamma,\omega)
\end{equation}
for any graph $\Gamma$ with orientation $\omega$. The differential $\partial$ on $\gc$ is given by the following formula:
\[ \partial(\Gamma):=\sum_{e \in E(\Gamma)} \Gamma/e \]
where $\Gamma$ is any oriented graph and the sum is taken over all edges of $\Gamma$. Note that relation \eqref{eqn_graphdummy} implies that any graph which contains a loop is zero in $\gc$ since there is an orientation reversing automorphism of such a graph given by permuting the half-edges of the loop; consequently there is no problem in the definition of $\partial$ which might otherwise arise by trying to contract edges that are loops.

There is a natural bigrading on $\gc$ where $\gc[ij]$ is the vector space generated by isomorphism classes of oriented graphs with $i$ vertices and $j$ edges. In this bigrading the differential has bidegree $(-1,-1)$.
\end{defi}

\begin{lemma}
The map $\partial:\gc[ij] \to \gc[i-1,j-1]$ is a differential, i.e. $\partial^2=0$.
\end{lemma}

\begin{proof}
This follows from the fact that contracting edges in the opposite order in an oriented graph produces opposite orientations on the resulting graph. Let $\Gamma$ be an oriented graph and $e$, $e'$ be two of its edges; say that $e\sim e'$ if and only if they have the same set of incident vertices, then
\begin{align*}
\partial^2(\Gamma) & = \sum_{\begin{subarray}{c} e,e' \in E(\Gamma): \\ e \nsim e' \\ \end{subarray}} (\Gamma/e)/e', \\
& = -\sum_{\begin{subarray}{c} e,e' \in E(\Gamma): \\ e \nsim e' \\ \end{subarray}} (\Gamma/e')/e = -\partial^2(\Gamma).
\end{align*}
It follows that $\partial^2(\Gamma)=0$.
\end{proof}

\begin{rem}
It was necessary to introduce the relation $\sim$ in order to omit the possibility of contracting edges which are loops, which is forbidden.
\end{rem}

This means that the object defined in Definition \ref{def_graph} is indeed a complex as claimed and we can therefore consider its homology which we call \emph{graph homology} and denote by $\gh$.

\section{Lie superalgebra homology} \label{sec_liehom}

In this section we will recall the basic definition of Lie superalgebra homology (cf. \cite{fuchscohom}) and introduce a series of Lie superalgebras $\galg{2n|m}$ which we intend to compute the (relative) stable Lie superalgebra homology of. This series of Lie superalgebras is the super-analogue of the series of Lie algebras introduced by Kontsevich in \cite{kontsg}. Informally, $\galg{2n|m}$ is the Lie superalgebra of Hamiltonian vector fields on a supermanifold of dimension $2n|m$.

\begin{defi}
We define $\mathfrak{h}_{2n|m}$ as the Poisson algebra of Hamiltonians on the supermanifold $\mathbb{C}^{2n|m}$, that is to say that:
\begin{enumerate}
\item
the underlying algebra of $\mathfrak{h}_{2n|m}$ is the supercommutative algebra
\[ S(\mathbb{C}^{2n|m}) = \mathbb{C}\langle p_1,\ldots,p_n;q_1,\ldots,q_n\rangle \otimes \Lambda(\Pi x_1,\ldots,\Pi x_m), \]
\item
$\mathfrak{h}_{2n|m}$ is equipped with the canonical Poisson bracket
\begin{equation} \label{eqn_poibra}
\{a,b\}:= \sum_{i=1}^n \left[\frac{\partial a}{\partial p_i}\frac{\partial b}{\partial q_i} - \frac{\partial a}{\partial q_i}\frac{\partial b}{\partial p_i}\right] - (-1)^{|a|}\sum_{i=1}^m\left[ \frac{\partial a}{\partial x_i}\frac{\partial b}{\partial x_i} \right]; \text{ for } a,b \in \mathfrak{h}_{2n|m}.
\end{equation}
\end{enumerate}

The Lie superalgebra $\galg{2n|m}$ is defined to be the Lie subalgebra of $\mathfrak{h}_{2n|m}$ consisting of polynomials whose terms have quadratic order or higher, that is to say polynomials with vanishing constant and linear terms. The Lie superalgebra $\tgalg{2n|m}$ is defined to be the Lie subalgebra of $\galg{2n|m}$ consisting of polynomials whose terms have cubic order or higher, that is to say polynomials with vanishing constant, linear and quadratic terms.
\end{defi}

\begin{rem} \label{rem_symsub}
The vector space $\mathbb{C}^{2n|m}$ sits inside the Lie superalgebra $\mathfrak{h}_{2n|m}$ as its linear part. The restriction of the Poisson bracket $\{-,-\}$ on $\mathfrak{h}_{2n|m}$ to the vector space $\mathbb{C}^{2n|m}$ is the canonical symplectic form $\langle-,-\rangle$ defined by equation \eqref{eqn_cansymfrm}. Consider the Lie subalgebra $\mathfrak{k}_{2n|m}$ of $\galg{2n|m}$ consisting of quadratic polynomials, i.e. those polynomials whose nonquadratic terms all vanish. The Lie superalgebra $\mathfrak{k}_{2n|m}$ acts on $\mathfrak{h}_{2n|m}$ via the adjoint action, in particular it acts on $\mathbb{C}^{2n|m}$, $\galg{2n|m}$ and $\tgalg{2n|m}$. It follows from the Jacobi identity that we have a map of Lie superalgebras
\begin{displaymath}
\begin{array}{ccc}
\mathfrak{k}_{2n|m} & \to & \osp{2n|m}, \\
a & \mapsto & [\{a,-\}: x \mapsto \{a,x\}, \text{ for } x \in \mathbb{C}^{2n|m}]. \\
\end{array}
\end{displaymath}
In fact this map is an isomorphism; this follows from the Poincar\'e Lemma \cite{ginzsg} after identifying the above map with the natural map from the space of Hamiltonians to the space of symplectic vector fields. The consequence of this is that $\galg{2n|m}$ splits as a semidirect product of $\osp{2n|m}$ and $\tgalg{2n|m}$;
\[ \galg{2n|m} = \osp{2n|m}\ltimes \tgalg{2n|m}, \]
where $\osp{2n|m}$ acts canonically on $\tgalg{2n|m}$ according to the Leibniz rule.
\end{rem}

We need to recall how the Chevalley-Eilenberg complex is defined for Lie superalgebras:

\begin{defi}
Let $\mathfrak{g}$ be a Lie superalgebra: the underlying space of the Chevalley-Eilenberg complex of $\mathfrak {g}$ is the exterior superalgebra $\Lambda_{\bullet}(\mathfrak{g})$ which is defined to be the quotient of $T(\mathfrak{g})$ by the ideal generated by the relation
\[ g \otimes h = -(-1)^{|g||h|} h \otimes g; \quad g,h \in \mathfrak{g}. \]
There is a natural grading on $\Lambda_{\bullet}(\mathfrak{g})$ where an element $g\in\mathfrak{g}$ has bidegree $(1,|g|)$ in $\Lambda_{\bullet\bullet}(\mathfrak{g})$ and total degree $|g|+1$; this implicitly defines a bigrading and total grading on the whole of $\Lambda_{\bullet}(\mathfrak{g})$.

The differential $d:\Lambda_{i}(\mathfrak{g}) \to \Lambda_{i-1}(\mathfrak{g})$ is defined by the following formula:
\[ d(g_1\wedge\ldots\wedge g_m):= \sum_{1\leq i < j \leq m} (-1)^{p(g)} [g_i,g_j]\wedge g_1\wedge\ldots\wedge \hat{g_i}\wedge\ldots\wedge \hat{g_j}\wedge\ldots\wedge g_m, \]
for $g_1,\ldots,g_m \in \mathfrak{g}$; where
\[ p(g):=|g_i|(|g_1|+\ldots+|g_{i-1}|)+|g_j|(|g_1|+\ldots+|g_{j-1}|)+|g_i||g_j|+i+j-1. \]
The differential $d$ has bidegree $(-1,0)$ in the above bigrading. We will denote the Chevalley-Eilenberg complex of $\mathfrak{g}$ by $C_{\bullet}(\mathfrak{g})$; the homology of $C_{\bullet}(\mathfrak{g})$ will be called the Lie superalgebra homology of $\mathfrak{g}$ and will be denoted by $H_{\bullet}(\mathfrak{g})$.
\end{defi}

\begin{rem}
Let $\mathfrak{g}$ be a Lie superalgebra; $\mathfrak{g}$ acts on $\Lambda(\mathfrak{g})$ via the adjoint action. This action commutes with the Chevalley-Eilenberg differential $d$; in fact it is nullhomotopic (cf. \cite{fuchscohom} or \cite{loday}).

As a consequence of this the space
\[ C_{\bullet}(\tgalg{2n|m})_{\osp{2n|m}} \]
of $\osp{2n|m}$-coinvariants of the Chevalley-Eilenberg complex of $\tgalg{2n|m}$ forms a complex when equipped with the Chevalley-Eilenberg differential $d$. This is the \emph{relative} Chevalley-Eilenberg complex of $\galg{2n|m}$ modulo $\osp{2n|m}$ (cf. \cite{fuchscohom}) and is denoted by $C_{\bullet}(\galg{2n|m},\osp{2n|m})$. The homology of this complex is called the \emph{relative} homology of $\galg{2n|m}$ modulo $\osp{2n|m}$ and is denoted by $H_{\bullet}(\galg{2n|m},\osp{2n|m})$.
\end{rem}

\begin{rem} \label{rem_ceordgrd}
There is a further grading on $C_{\bullet}(\galg{2n|m},\osp{2n|m})$ determined by the grading by order on $\galg{2n|m}$. In this grading an element $g \in \galg{2n|m}$ has bidegree $(1,o)$ where $o$ is the order of the homogeneous element $g$. The differential $d$ has bidegree $(-1,-2)$ in this grading.
\end{rem}

\begin{rem}
Consider the Lie algebra $\galg{2n|0} = \mathfrak{sp}_{2n}\ltimes \tgalg{2n|0}$. By considering the Hochschild-Serre spectral sequence associated to the Lie subalgebra $\mathfrak{sp}_{2n}$ of $\galg{2n|0}$, it is possible to express the Lie algebra homology of $\galg{2n|0}$ as a tensor product of the Lie algebra homology of $\mathfrak{sp}_{2n}$ and the relative Lie algebra homology of $\galg{2n|0}$ modulo $\mathfrak{sp}_{2n}$. This follows from standard results in the representation theory of Lie algebras; namely that $\mathfrak{sp}_{2n}$ is a simple Lie algebra and that the Lie algebra homology of a simple Lie algebra with coefficients in a simple nontrivial module over that Lie algebra vanishes, cf. \cite[\S 7.8]{Weibel}. This calculation was originally outlined in \cite{kontsg}, except that the spectral sequence calculations were performed via a filtration on the graph complex using ad hoc arguments.

If $m > 0$ then it may be that the Lie superalgebra homology of $\galg{2n|m}$ does not decompose in this way; this is because the representation theory of Lie \emph{superalgebras} (cf. \cite{schlie}) is substantially different from the representation theory of Lie \emph{algebras} which makes it difficult to apply standard arguments from the ungraded theory to Lie superalgebras.
\end{rem}

\section{Invariant theory for the Lie superalgebras $\osp{2n|m}$} \label{sec_invthy}

In the first half of this section we will summarise the relevant results on the invariant theory for the Lie superalgebras $\osp{2n|m}$ that we will need to make use of in the next section. Here we will be very brief and only list the results on the invariant theory of $\osp{2n|m}$ that are required for our constructions in section \ref{sec_gphhommain}, rather than discussing the invariant theory of $\osp{2n|m}$ in all its detail. The main reference for the invariant theory for $\osp{2n|m}$ is \cite{sergeev}. It will also be necessary to translate these results to the dual framework of coinvariants and this is carried out in the second half of the section.

Let $V:=\mathbb{C}^{2n|m}$ and let $\kappa:V^{\otimes 2k} \to \mathbb{C}$ be the map defined by the formula
\[ \kappa(x_1\otimes\ldots\otimes x_{2k}):= \langle x_1,x_2 \rangle\langle x_3,x_4 \rangle\ldots\langle x_{2k-1},x_{2k}\rangle, \]
where $\langle -,-\rangle$ is the canonical symplectic form \eqref{eqn_cansymfrm}. It follows from the definitions that $\kappa$ is $\osp{2n|m}$-invariant. The following definition produces the fundamental invariants of the Lie superalgebra $\osp{2n|m}$:

\begin{defi} \label{def_invdef}
Let $k$ be a positive integer and $c \in \mathscr{OC}(k)$ be an oriented chord diagram. Write $c$ as
\begin{equation} \label{eqn_invdefdummy}
c:= (i_1,j_1),\ldots,(i_k,j_k)
\end{equation}
We define a permutation $\sigma_c \in S_{2k}$ by
\begin{equation} \label{eqn_invdefdummya}
\begin{array}{cccccccc}
& i_1 & j_1 & i_2 & j_2 & & i_k & j_k \\
\sigma_c:= & \downarrow & \downarrow & \downarrow & \downarrow & \cdots & \downarrow & \downarrow \\
& 1 & 2 & 3 & 4 & & 2k-1 & 2k \\
\end{array};
\end{equation}
this permutation is not uniquely defined for any given chord diagram $c$ and depends upon the choice that was made when writing $c$ down in the form \eqref{eqn_invdefdummy}, namely an ordering of the pairs in the chord diagram.

We define the map $\beta_c:V^{\otimes 2k} \to \mathbb{C}$ by the formula
\[ \beta_c(x):= \kappa(\sigma_c\cdot x). \]
Although $\sigma_c$ will depend upon the choice we made when we wrote $c$ down in \eqref{eqn_invdefdummy}, $\beta_c$ will be independent of this choice since the symplectic form $\langle -,- \rangle$ is even, hence $\beta_c$ is well defined.
\end{defi}

\begin{rem} \label{rem_chdinv}
Note that $\beta_c$ is $\osp{2n|m}$-invariant as for all $\xi \in \osp{2n|m}$ and $x\in V^{\otimes 2k}$,
\[ \beta_c(\xi\cdot x) = \kappa(\sigma_c\cdot\xi\cdot x) = \kappa(\xi \cdot \sigma_c \cdot x) = 0. \]

Also note that for all $\sigma \in S_{2k}$,
\[ \beta_{\sigma\cdot c}(\sigma\cdot x) = \beta_c(x). \]
If $\tau:=(i_r,j_r)\in S_{2k}$, $1\leq r \leq k$; then
\[ \beta_c(\tau\cdot x) = -\beta_c(x), \]
as the symplectic form $\langle -,- \rangle$ is antisymmetric.
\end{rem}

Given a chord diagram $c \in \mathscr{C}(k)$ written as
\[ c:= \{i_1,j_1\},\ldots,\{i_k,j_k\}, \]
we can canonically define an \emph{oriented} chord diagram $\hat{c}$ as follows: we may assume that $i_r < j_r$ for all $r$, then $\hat{c}$ is defined as
\[ \hat{c}:= (i_1,j_1),\ldots,(i_k,j_k). \]
We can now state the fundamental results from the invariant theory of $\osp{2n|m}$, suitably translated from \cite{sergeev}, that we will need to make use of later on:

\begin{lemma} \label{lem_invthy}
Let $V:=\mathbb{C}^{2n|m}$:
\begin{enumerate}
\item
If $k\leq n$ then $[(V^*)^{\otimes 2k-1}]^{\osp{2n|m}} = 0$.
\item
If $k \leq n$ then $[(V^*)^{\otimes 2k}]^{\osp{2n|m}}$ has a basis indexed by chord diagrams:
\[ \beta_{\hat{c}}, \ c \in \mathscr{C}(k). \]
Moreover, if $k\leq n$ then $OSp_{2n|m}$ acts trivially on $[(V^*)^{\otimes 2k}]^{\osp{2n|m}}$.
\end{enumerate}
\end{lemma}
\noproof

The rest of this section will describe how to translate the above result to the dual framework of coinvariants. We begin with a definition of the fundamental coinvariants of $\osp{2n|m}$:

\begin{defi} \label{def_coinv}
Let $V:= \mathbb{C}^{2n|m}$, let $k\leq n$ be a positive integer and let $c \in \mathscr{OC}(k)$ be an oriented chord diagram: we define an element $u^{2n|m}_k(c) \in (V^{\otimes 2k})_{\osp{2n|m}}$ by the formula
\[ u^{2n|m}_k(c):= \sigma_c^{-1} \cdot (p_1\otimes q_1)\otimes\ldots\otimes(p_k\otimes q_k), \]
where $\sigma_c$ is given by equation \eqref{eqn_invdefdummya}.
\end{defi}

\begin{rem}
It follows from Lemma \ref{lem_invthy} part (ii) that $u^{2n|m}_k(c)$ is well defined; given an oriented chord diagram $c \in \mathscr{OC}(k)$, the permutation $\sigma_c$ is not well defined. However, since
\begin{equation} \label{eqn_invcovdual}
\left[(V^{\otimes 2k})_{\osp{2n|m}}\right]^* = \left[(V^*)^{\otimes 2k}\right]^{\osp{2n|m}}
\end{equation}
it follows from Lemma \ref{lem_invthy} that $OSp_{2n|m}$ acts trivially on $(V^{\otimes 2k})_{\osp{2n|m}}$; this in turn implies the identity:
\[ (p_1\otimes q_1)\otimes\ldots\otimes(p_k\otimes q_k) = (p_{\tau(1)}\otimes q_{\tau(1)})\otimes\ldots\otimes(p_{\tau(k)}\otimes q_{\tau(k)}) \mod \left[\osp{2n|m},V^{\otimes 2k}\right], \]
for all $\tau \in S_n$. It follows from this identity that $u^{2n|m}_k(c)$ is independent of the choice made in defining $\sigma_c$.
\end{rem}

We are now ready to give a description for the coinvariants of $\osp{2n|m}$ which is dual to the description of invariants formulated in Lemma \ref{lem_invthy}:

\begin{lemma} \label{lem_covbas}
Let $V:=\mathbb{C}^{2n|m}$ and let $k \leq n$ be a positive integer, then $(V^{\otimes 2k})_{\osp{2n|m}}$ has a basis indexed by chord diagrams:
\[ u^{2n|m}_k(\hat{c}), \ c \in \mathscr{C}(k). \]
Moreover, this basis is dual to the basis of $[(V^*)^{\otimes 2k}]^{\osp{2n|m}}$ given by
\[ \beta_{\hat{c}}, \ c \in \mathscr{C}(k). \]
\end{lemma}

\begin{proof}
Given a chord diagram $c \in \mathscr{C}(k)$, we see that
\begin{displaymath}
\begin{split}
\beta_{\hat{c}}\left[u^{2n|m}_k(\hat{c})\right] & = \beta_{\hat{c}}\left[\sigma_{\hat{c}}^{-1} \cdot (p_1\otimes q_1)\otimes\ldots\otimes(p_k\otimes q_k)\right] \\
& = \kappa\left[(p_1\otimes q_1)\otimes\ldots\otimes(p_k\otimes q_k)\right]=1.
\end{split}
\end{displaymath}

Let $c' \in \mathscr{C}(k)$ be another chord diagram such that $c'\neq c$ and write
\begin{displaymath}
\begin{split}
& \hat{c} = (i_1,j_1)\ldots(i_k,j_k), \\
& \hat{c'} = (i'_1,j'_1)\ldots(i'_k,j'_k); \\
\end{split}
\end{displaymath}
since $c'\neq c$ we may assume that $i_1 \in \{i'_1,j'_1\}$ and $j_1 \in \{i'_2,j'_2\}$. It follows that
\begin{displaymath}
\begin{split}
\beta_{\hat{c'}}\left[u^{2n|m}_k(\hat{c})\right] & = \kappa\left[\sigma_{\hat{c'}}\sigma_{\hat{c}}^{-1}\cdot (p_1\otimes q_1)\otimes\ldots\otimes(p_k\otimes q_k) \right] \\
& = \pm \langle p_1,x \rangle\times\ldots \ ; \ x=p_j \text{ or } x=q_j \text{ for some } j\geq 2 \\
& = 0. \\
\end{split}
\end{displaymath}
Now it follows from equation \eqref{eqn_invcovdual} and Lemma \ref{lem_invthy} that $u^{2n|m}_k(\hat{c})$, $c \in \mathscr{C}(k)$ forms a basis of $(V^{\otimes 2k})_{\osp{2n|m}}$.
\end{proof}

\section{The main theorem: an identification of graph homology as the Lie superalgebra homology of supersymplectic vector fields} \label{sec_gphhommain}

In this section we shall prove our main theorem: that the homology of the graph complex $\gc$ can be calculated as the homology of the complex
\[ \mathcal{C}_{\bullet}:= \dilim{n,m}{C_{\bullet}(\galg{2n|m},\osp{2n|m})}; \]
in fact we shall actually show that these two complexes are isomorphic. This is the super-analogue of Kontsevich's well-known theorem on graph homology \cite{kontsg}. The main tool that we will need to make use of is the invariant theory of the Lie superalgebras $\osp{2n|m}$ which we established in section \ref{sec_invthy}. Using the $\osp{2n|m}$-invariants described in the preceding section we will construct a map between the complexes $\mathcal{C}_{\bullet}$ and $\gc$. The results from section \ref{sec_invthy} on the invariant theory of $\osp{2n|m}$ will then be applied to show that this map is an isomorphism.

The first point to mention regards the grading on $\mathcal{C}_\bullet$. As a corollary of Lemma \ref{lem_invthy} and equation \eqref{eqn_invcovdual} we see that all those elements of \emph{odd} order (cf. Remark \ref{rem_ceordgrd}) in the complex
\[ \mathcal{C}_{\bullet}:= \dilim{n,m}{C_{\bullet}(\galg{2n|m},\osp{2n|m})}. \]
vanish in the limit. This gives us a bigrading on $\mathcal{C}$ given by
\[ \mathcal{C}_{ij}:= \dilim{n,m}{C_{i,2j}(\galg{2n|m},\osp{2n|m})}. \]
where $C_{i,2j}(\galg{2n|m},\osp{2n|m})$ consists of all those elements in $\left[\Lambda^i(\tgalg{2n|m})\right]_{\osp{2n|m}}$ which have order $2j$. The differential $d$ has bidegree $(-1,-1)$ in this grading.

Now we will go about constructing an isomorphism between $\mathcal{C}_{\bullet\bullet}$ and $\gc$. The first step is to describe how to construct an oriented graph from an oriented chord diagram:

\begin{defi} \label{def_grpchd}
Let $k_1,\ldots,k_i$ be positive integers such that $k_1+\ldots+k_i = 2k$ and let
\[ c:=(i_1,j_1),\ldots,(i_k,j_k)\in \mathscr{OC}(k) \]
be an oriented chord diagram; we define an oriented graph $\Gamma_{k_1,\ldots,k_i}(c)$ with half edges $h_1,\ldots,h_{2k}$ as follows:
\begin{enumerate}
\item
The vertices of $\Gamma$ are
\begin{equation} \label{eqn_grpchddummy}
\{h_1,\ldots,h_{k_1}\},\{h_{k_1+1},\ldots,h_{k_1+k_2}\},\ldots,\{h_{k_1+\ldots+k_{i-1}+1},\ldots,h_{k_1+\ldots+k_i}\}
\end{equation}
\item
The edges of $\Gamma$ are
\begin{equation} \label{eqn_grpchddummya}
(h_{i_1},h_{j_1}),\ldots,(h_{i_k},h_{j_k})
\end{equation}
\item
The orientation on $\Gamma$ is given by ordering the vertices as in \eqref{eqn_grpchddummy} and orienting the edges as in \eqref{eqn_grpchddummya}.
\end{enumerate}
\end{defi}

This construction satisfies some simple identities which are easily verified:

\begin{lemma} \label{lem_grpchdid}
Let $k_1,\ldots,k_i$ be positive integers such that $k_1+\ldots+k_i = 2k$ and let
\[ c:=(i_1,j_1),\ldots,(i_k,j_k)\in \mathscr{OC}(k) \]
be an oriented chord diagram:
\begin{enumerate}
\item
Let $\tau:=(i_r,j_r) \in S_{2k}$, $1\leq r \leq k$; then
\[ \Gamma_{k_1,\ldots,k_i}(\tau\cdot c) = - \Gamma_{k_1,\ldots,k_i}(c). \]
\item
Let $\sigma:=(\sigma_1,\ldots,\sigma_i) \in S_{k_1}\times\ldots\times S_{k_i} \subset S_{2k}$, then
\[ \Gamma_{k_1,\ldots,k_i}(\sigma\cdot c) = \Gamma_{k_1,\ldots,k_i}(c). \]
\item
Let $\sigma \in S_i$, then
\[ \Gamma_{k_{\sigma(1)},\ldots,k_{\sigma(i)}}((\sigma_{k_1,\ldots,k_i})^{-1}\cdot c) = \sgn(\sigma) \Gamma_{k_1,\ldots,k_i}(c). \]
\end{enumerate}
\end{lemma}
\noproof

Now we can use this construction and the $\osp{2n|m}$-invariants defined in Definition \ref{def_invdef} to define a map between $C_{\bullet\bullet}(\galg{2n|m},\osp{2n|m})$ and $\gc$:

\begin{defi} \label{def_liegphmap}
\
\begin{enumerate}
\item
The map $\Phi_{ij}^{2n|m}:C_{ij}(\galg{2n|m},\osp{2n|m}) \to \gc[ij]$ is defined by the following formula: let
\[ x:=(x_{11}\ldots x_{1k_1})\wedge(x_{21}\ldots x_{2k_2})\wedge\ldots\wedge(x_{i1}\ldots x_{ik_i}) \in C_{ij}(\galg{2n|m},\osp{2n|m}), \]
where $x_{pq} \in \mathbb{C}^{2n|m}$; then
\begin{equation} \label{eqn_liegphmap}
\Phi_{ij}^{2n|m}(x):=\sum_{c\in\mathscr{C}(j)}\beta_{\hat{c}}(x)\Gamma_{k_1,\ldots,k_i}(\hat{c}).
\end{equation}
\item \label{def_liegphmap_limmap}
The map $\Phi:\mathcal{C}_{\bullet\bullet} \to \gc$ is given by the formula
\[ \Phi:= \dilim{n,m}{\Phi^{2n|m}}. \]
\end{enumerate}
\end{defi}

\begin{rem}
It follows from Lemma \ref{lem_grpchdid} and Remark \ref{rem_chdinv} that $\Phi$ is well defined. For instance let
\[ \sigma:=(\sigma_1,\ldots,\sigma_i) \in S_{k_1}\times\ldots\times S_{k_i} \subset S_{2k}; \]
for any chord diagram $c:=\{i_1,j_1\},\ldots,\{i_k,j_k\}\in \mathscr{C}(k)$ where $i_r < j_r$ for all $r$, let $N^\sigma_c$ be the number of integers $r$ between $1$ and $k$ such that $\sigma(i_r) > \sigma(j_r)$; then
\begin{displaymath}
\begin{split}
\Phi(\sigma\cdot x) & = \sum_{c\in\mathscr{C}(k)}\beta_{\hat{c}}(\sigma\cdot x)\Gamma_{k_1,\ldots,k_i}(\hat{c}) \\
& = \sum_{c\in\mathscr{C}(k)}\beta_{\widehat{\sigma\cdot c}}(\sigma\cdot x)\Gamma_{k_1,\ldots,k_i}(\widehat{\sigma\cdot c}) \\
& = \sum_{c\in\mathscr{C}(k)}(-1)^{N^\sigma_c}\beta_{\hat{c}}(x)(-1)^{N^\sigma_c}\Gamma_{k_1,\ldots,k_i}(\hat{c}) = \Phi(x). \\
\end{split}
\end{displaymath}
This shows that $\Phi$ is well defined modulo the symmetry relations in $\tgalg{2n|m}$. A similar calculation shows that it is well defined modulo the antisymmetry relations present in $\Lambda(\tgalg{2n|m})$. It follows from Remark \ref{rem_chdinv} that it is $\osp{2n|m}$-invariant.
\end{rem}

Next we need to show that our map $\Phi$ respects the differentials in $C_{\bullet\bullet}(\galg{2n|m},\osp{2n|m})$ and $\gc$:

\begin{prop} \label{prop_mapcpx}
The map $\Phi^{2n|m}:C_{\bullet\bullet}(\galg{2n|m},\osp{2n|m}) \to \gc$ is a map of complexes, that is to say that
\[ \Phi^{2n|m} \circ d = \partial \circ \Phi^{2n|m}. \]
\end{prop}

\begin{proof}
Let
\[ x:=y_1\wedge\ldots\wedge y_p \in \Lambda^p(\tgalg{2n|m}), \quad y_r:=x_{N_r+1}\ldots x_{N_r+k_r}; \]
where $k_1,\ldots,k_p$ are positive integers such that $k_1+\ldots+k_p=2k$, $N_t:=k_1+\ldots+k_{t-1}$ and $x_r \in \mathbb{C}^{2n|m}$. A straightforward application of the Leibniz rule shows that
\[ \{ y_i,y_j \} = \sum_{s=1}^{k_i}\sum_{t=1}^{k_j} (-1)^{\tilde{q}} \langle x_{N_i+s},x_{N_j+t} \rangle x_{N_i+1}\ldots \widehat{x_{N_i+s}}\ldots x_{N_i+k_i}x_{N_j+1}\ldots \widehat{x_{N_j+t}} \ldots x_{N_j+k_j}, \]
where $\tilde{q}:=|x_{N_j+s}|(|x_{N_i+s+1}|+\ldots+|x_{N_i+k_i}|+|x_{N_j+1}|+\ldots+|x_{N_j+t-1}|)$. It follows that
\begin{displaymath}
d(x) = \sum_{\begin{subarray}{c} 1\leq i < j \leq p \\ 1\leq s \leq k_i \\ 1 \leq t \leq k_j \end{subarray}} \left[\begin{array}{l} (-1)^q \langle x_{N_i+s},x_{N_j+t} \rangle\left(x_{N_i+1}\ldots \widehat{x_{N_i+s}}\ldots x_{N_i+k_i}x_{N_j+1}\ldots \widehat{x_{N_j+t}} \ldots x_{N_j+k_j}\right)\wedge \\ y_1\wedge\ldots\wedge\hat{y_i}\wedge\ldots\wedge\hat{y_j}\wedge\ldots\wedge y_p \end{array}\right],
\end{displaymath}
where $q:=|y_i|(|y_1|+\ldots+|y_{i-1}|)+ |y_j|(|y_1|+\ldots+|y_{j-1}|)+|y_i||y_j|+i+j-1+\tilde{q}$.

This gives us the following identity:
\begin{equation} \label{eqn_mapcpxdummy}
\Phi d(x) = \sum_{\begin{subarray}{c} 1\leq i < j \leq p \\ 1\leq s \leq k_i \\ 1 \leq t \leq k_j \end{subarray}} \: \sum_{c\in\mathscr{C}(k-1)} \left[(-1)^q \langle x_{N_i+s},x_{N_j+t}\rangle \beta_{\hat{c}}[x(i,j;s,t)] \Gamma_{k_i+k_j-2,k_1,\ldots,\hat{k_i},\ldots,\hat{k_j},\ldots,k_p}(\hat{c})\right],
\end{equation}
where
\[ x(i,j;s,t):=\left(x_{N_i+1}\ldots \widehat{x_{N_i+s}}\ldots x_{N_i+k_i}x_{N_j+1}\ldots \widehat{x_{N_j+t}} \ldots x_{N_j+k_j}\right)\otimes y_1\otimes\ldots\otimes\hat{y_i}\otimes\ldots\otimes\hat{y_j}\otimes\ldots\otimes y_p. \]

Let $c\in\mathscr{C}(k)$ be a chord diagram, then
\[ \partial[\Gamma_{k_1,\ldots k_p}(\hat{c})] = \sum_{\begin{subarray}{c} 1\leq i < j \leq p \\ 1\leq s \leq k_i \\ 1\leq t \leq k_j, \\ \text{such that} \\ \{N_i+s,N_j+t\}\in c \end{subarray}} \Gamma_{k_1,\ldots k_p}(\hat{c})/(h_{N_i+s},h_{N_j+t}). \]
It follows that
\begin{equation} \label{eqn_mapcpxdummya}
\partial\Phi(x) = \sum_{\begin{subarray}{c} 1\leq i < j \leq p \\ 1\leq s \leq k_i \\ 1\leq t \leq k_j \end{subarray}}\sum_{\begin{subarray}{c} c \in \mathscr{C}(k): \\ \{N_i+s,N_j+t\}\in c \\ \end{subarray}} \beta_{\hat{c}}(x)[\Gamma_{k_1,\ldots k_p}(\hat{c})/(h_{N_i+s},h_{N_j+t})].
\end{equation}

By comparing equation \eqref{eqn_mapcpxdummy} to equation \eqref{eqn_mapcpxdummya} it is easy to verify that
\[ \Phi d(x) = \partial\Phi(x). \]
\end{proof}

Having constructed the map $\Phi:\mathcal{C}_{\bullet\bullet} \to \gc$ we now need to show that $\Phi$ is an isomorphism. This is done by constructing an inverse to $\Phi$ using the description of $\osp{2n|m}$-coinvariants we derived in the last section, cf. Definition \ref{def_coinv} and Lemma \ref{lem_covbas}.

Let $V:=\mathbb{C}^{2n|m}$; given a vector $u \in V^{\otimes 2k}$ and positive integers $k_1,\ldots,k_i \geq 3$ such that $k_1+\ldots+k_i=2k$, let $\lceil u\rceil_{k_1,\ldots,k_i}$ denote the canonical image of $u$ under the map
\[ V^{\otimes 2k} \to S^{k_1}(V)\otimes\ldots\otimes S^{k_i}(V) \to \Lambda^i(\tgalg{2n|m}). \]
We are now ready to construct the inverse to the map $\Phi$:

\begin{defi}
The map $\Psi^{2n|m}_{ij}: \mathcal{G}_{ij} \to C_{ij}(\galg{2n|m},\osp{2n|m})$ is defined for $j\leq n$ as follows: given a graph $\Gamma \in \mathcal{G}_{ij}$, let $c \in \mathscr{OC}(j)$ be an oriented chord diagram and $k_1,\ldots,k_i$ be positive integers such that $\Gamma = \Gamma_{k_1,\ldots,k_i}(c)$\footnote{Such a chord diagram and sequence of positive integers always exist.}; then
\[ \Psi^{2n|m}_{ij}(\Gamma):=\lceil u^{2n|m}_j(c)\rceil_{k_1,\ldots,k_i}. \]
\end{defi}

\begin{rem}
Firstly, it is worth mentioning that $\Psi^{2n|m}_{ij}$ is defined as a map of vector spaces and not as a map of complexes. Secondly $\Psi^{2n|m}_{ij}(\Gamma)$ is independent of the choice of the oriented chord diagram $c$ and sequence of positive integers $k_1,\ldots,k_i$; for instance, since
\[ p_i \otimes q_i = - q_i \otimes p_i \mod OSp_{2n|m}, \]
it follows from Lemma \ref{lem_invthy} that different choices for the orientation on the pairs in the chord diagram $c$ do not affect the value of $\Psi^{2n|m}_{ij}(\Gamma)$. Furthermore, the symmetry and antisymmetry relations present in $\tgalg{2n|m}$ and $\Lambda^i(\tgalg{2n|m})$ guarantee that different choices of the chord diagram $c$ and sequence of positive integers $k_1,\ldots,k_i$ do not affect the value of $\Psi^{2n|m}_{ij}(\Gamma)$.
\end{rem}

\begin{lemma} \label{lem_gphiso}
The following identities hold for $j\leq n$:
\begin{align}
\label{eqn_gphisodummy}
\Psi^{2n|m}_{ij} \circ \Phi^{2n|m}_{ij} & = \id_{C_{ij}(\galg{2n|m},\osp{2n|m})}, \\
\label{eqn_gphisodummya}
\Phi^{2n|m}_{ij} \circ \Psi^{2n|m}_{ij} & = \id_{\mathcal{G}_{ij}}.
\end{align}
\end{lemma}

\begin{proof}
Let $c \in \mathscr{C}(j)$ be a chord diagram and $k_1,\ldots,k_i$ be positive integers such that $k_1+\ldots+k_i=2j$,
\begin{displaymath}
\begin{split}
\Psi^{2n|m}_{ij} \Phi^{2n|m}_{ij}\left(\lceil u^{2n|m}_j(\hat{c})\rceil_{k_1,\ldots,k_i}\right) & = \Psi^{2n|m}_{ij}\left(\sum_{c'\in\mathscr{C}(j)}\beta_{\hat{c'}}\left[u^{2n|m}_j(\hat{c})\right] \Gamma_{k_1\ldots,k_i}(\hat{c'})\right) \\
& = \Psi^{2n|m}_{ij}\left(\Gamma_{k_1,\ldots,k_i}(\hat{c})\right) \\
& = \lceil u^{2n|m}_j(\hat{c})\rceil_{k_1,\ldots,k_i}. \\
\end{split}
\end{displaymath}
Equation \eqref{eqn_gphisodummy} now follows from Lemma \ref{lem_covbas}.

Again, let $c \in \mathscr{C}(j)$ be a chord diagram and $k_1,\ldots,k_i$ be positive integers such that $k_1+\ldots+k_i=2j$,
\begin{displaymath}
\begin{split}
\Phi^{2n|m}_{ij} \Psi^{2n|m}_{ij}\left(\Gamma_{k_1,\ldots,k_i}(\hat{c})\right) & = \Phi^{2n|m}_{ij}\left(\lceil u^{2n|m}_j(\hat{c})\rceil_{k_1,\ldots,k_i}\right) \\
& = \sum_{c' \in \mathscr{C}(j)} \beta_{\hat{c'}}\left[u^{2n|m}_j(\hat{c})\right] \Gamma_{k_1,\ldots,k_i}(\hat{c'}) \\
& = \Gamma_{k_1,\ldots,k_i}(\hat{c}). \\
\end{split}
\end{displaymath}
This establishes equation \eqref{eqn_gphisodummya}.
\end{proof}

We are now in a position to prove the main theorem of this paper which is the super-analogue of Kontsevich's theorem on graph homology \cite{kontsg}:

\begin{theorem} \label{thm_gphiso}
The map $\Phi:\mathcal{C}_{\bullet\bullet} \to \gc$ defined in Definition \ref{def_liegphmap}(\ref{def_liegphmap_limmap}) is an isomorphism.
\end{theorem}

\begin{proof}
Recall that $\mathcal{C}_{\bullet\bullet}$ was defined as
\[ \mathcal{C}_{\bullet\bullet}:= \dilim{n,m}{C_{\bullet\bullet}(\galg{2n|m},\osp{2n|m})}. \]
It follows as a simple and formal consequence of Lemma \ref{lem_gphiso} that the map $\Phi$ is an isomorphism.
\end{proof}

\begin{cor}
\[ \gh = \dilim{n,m}{H_{\bullet\bullet}(\galg{2n|m},\osp{2n|m})}. \]
\end{cor}
\noproof

\begin{rem}
Given two graphs $\Gamma$ and $\Gamma'$, one can define their disjoint union $\Gamma \sqcup \Gamma'$ in an obvious way. If the graphs $\Gamma$ and $\Gamma'$ are oriented graphs then this induces a canonical orientation on $\Gamma \sqcup \Gamma'$. Since every graph is the disjoint union of its connected components, this gives the graph complex $\gc$ the structure of a commutative cocommutative differential graded Hopf algebra in which the multiplication is given by the disjoint union of graphs and the subcomplex of primitive elements is the space spanned by \emph{connected} graphs.

Furthermore, the complex
\[ \mathcal{C}_{\bullet\bullet}:= \dilim{n,m}{C_{\bullet\bullet}(\galg{2n|m},\osp{2n|m})} \]
has the canonical structure of a commutative cocommutative differential graded Hopf algebra; the comultiplication is given by the diagonal map
\begin{displaymath}
\begin{array}{ccc}
\tgalg{2n|m} & \to & \tgalg{2n|m} \oplus \tgalg{2n|m}, \\
g & \mapsto & (g,g); \\
\end{array}
\end{displaymath}
whilst the multiplication is induced by the canonical inclusion
\[ \tgalg{2n|m} \oplus \tgalg{2n'|m'} \hookrightarrow \tgalg{(2n+2n')|(m+m')}. \]

It is not hard to show that the map $\Phi:\mathcal{C}_{\bullet\bullet} \to \gc$ is a map of Hopf algebras, in particular the homology of the subcomplex of the graph complex which is spanned by \emph{connected} graphs may be calculated as the primitive homology of $\mathcal{C}_{\bullet\bullet}$. In order to keep this paper within reasonable limits we refrain from a more detailed discussion of the Hopf algebra structure; this matter will be dealt with more thoroughly in a subsequent paper in preparation by the author in collaboration with Andrey Lazarev.
\end{rem}

\end{document}